\documentclass[aap,MSNbibl,nameyear,dvips]{arximspdf}

\doi{10.1214/12-AAP858} %kopijuoti is PTS
\volume{23}
\issue{2}
\pubyear{2013}
\firstpage{834}
\lastpage{857}

\makeatletter

\newcommand{\eqref}[1]{(\ref{#1})}

\newcommand{\E}{\mathbb E}
\newcommand{\R}{\mathbb R}
\newcommand{\F}{\mathcal F}
\renewcommand{\P}{{\mathbb P}}
\newcommand{\Q}{{\mathbb Q}}

\newtheorem{theorem}{Theorem}[section]
\newtheorem{lemma}[theorem]{Lemma}
\newtheorem{corollary}[theorem]{Corollary}
\newtheorem{proposition}[theorem]{Proposition}

\newproclaim{condition}[theorem]{Condition}
\newproclaim{definition}[theorem]{Definition}

\makeatother

\begin{document}
\begin{frontmatter}

\title{Error distributions for random grid approximations of
multidimensional stochastic integrals\thanksref{T1}}
\runtitle{Error distributions}

\thankstext{T1}{Supported by the Swedish Foundation for
Strategic Research through the Gothenburg
Mathematical Modelling Centre.}

\begin{aug}
\author[A]{\fnms{Carl} \snm{Lindberg}\corref{}\ead[label=e1]{carl.lindberg@alumni.chalmers.se}}
\and
\author[A]{\fnms{Holger} \snm{Rootz\'en}\ead[label=e2]{hrootzen@chalmers.se}\ead[label=u1,url]{http://www.math.chalmers.se/\textasciitilde rootzen/}}
\runauthor{C. Lindberg and H. Rootz\'en}
\affiliation{Chalmers University of Technology and
Gothenburg University}
\address[A]{Department of Mathematical Sciences\\
Chalmers University of Technology\\
and \\
Gothenburg University\\
SE-412 96 G\"oteborg\\
Sweden\\
\printead{e1}\\
\hphantom{E-mail: }\printead*{e2}\\
\printead{u1}} %adresu isvedimo komanda gale!
\end{aug}

% HISTORY:
\received{\smonth{9} \syear{2011}}
\revised{\smonth{3} \syear{2012}}

% ABSTRACT
%
\begin{abstract}
This paper proves joint convergence of the approximation error for
several stochastic integrals with respect to local Brownian
semimartingales, for nonequidistant and random grids. The conditions
needed for convergence are that the Lebesgue integrals of the
integrands tend uniformly to zero and that the squared variation and
covariation processes converge. The paper also provides tools which
simplify checking these conditions and which extend the range for the
results. These results are used to prove an explicit limit theorem for
random grid approximations of integrals based on solutions of
multidimensional SDEs, and to find ways to ``design'' and optimize the
distribution of the approximation error. As examples we briefly
discuss strategies for discrete option hedging.
\end{abstract}

% KEYWORDS
%

\end{frontmatter}

%s1 #&#
\section{Introduction}\label{intro}
The error in numerical approximations of stochastic integrals is a
random variable, or, if one also is interested in the ``time''
development of the error, a stochastic process. Hence the most precise
evaluation of the error, which is possible to obtain, is to derive the
distribution of the error. The prototype example is the Euler method
for the stochastic integral $\int_0^t f(B(s),s) \,dB(s)$, for a
Brownian motion $B$. The Euler method approximates the integrand with a
step-function which is constant between the ``evaluation times'' (or,
in finance terminology, ``intervention times'') of the grid $i/n; i=0,
1, \ldots\,$. This leads to the approximation $\int_0^t f\circ\eta_n
\,dB(s)$, with $\eta_n(t)= i/n$ on the intervals $[i/n, (i+1)/n)$. In
\citet{Rootzen} it is shown that the approximation error\vadjust{\goodbreak} $U^n =
n^{1/2}\int_0^t(f - f\circ\eta_n)\,dB(s)$ converges stably in
distribution,
\[
U^n \Rightarrow_s \frac{1}{\sqrt{2}}\int_0^t f'(B(s),s)\,dW(s),
\]
where $W$ is a Brownian motion independent of $B$ and $f'(x,y) = \frac
{\partial f(x,y)}{\partial x}$, and where R\'enyi's quite useful
concept of stable convergence means that $U^n$ converges jointly with
any sequence which converges in probability.

The intuition behind this result is that ``the small wiggles of a
Brownian path are asymptotically independent of the global behavior of
the path.'' The result has seen much further development, in particular,
to the error in numerical solution schemes for SDEs, and has recently
found significant application in measuring the risks associated with
discrete hedging. A brief overview of some of this literature is given below.

The present paper generalizes this result in three ways: to joint
convergence of the approximation error for several stochastic
integrals, to local Brownian semimartingales instead of Brownian
motions, and to nonequidistant and random evaluation times. The tools
which help us quantify the intuition given above is \mbox{Girsanov's} theorem
which shows how a multidimensional Brownian motion is affected by a
change of measure, and L\'evy's characterization of a multidimensional
Brownian motion in terms of its square variation processes.

The conditions needed for convergence apply more generally than to
approximation schemes. They are that the Lebesgue integrals of the
integrands tend uniformly to zero in probability and that the square
variation and covariation processes converge in probability. We
additionally provide tools which simplify checking these conditions and
which extend the range of the results. Further we apply these results
to prove an explicit limit theorem for approximations of integrals
based on solutions of multidimensional SDEs.

One center of interest for this paper is the possibility to improve
approximation by using variable and random grids. In particular we
study approximation schemes where the evaluation times $i/n$ are
replaced by time points given by the recursion $\tau_0^n = 0$ and
\[
\tau_{k+1}^n = \tau_k^n + \frac{1}{n\theta(\tau_k^n)}
\]
for a positive adapted process $\theta(t)$.
We also study how the function $\theta$ can be chosen to design the
approximation error so that it has desirable properties. For example,
these could be homogeneous evolution of risk, or how to make the
approximation error have minimal standard deviation.

A main motivation for writing this paper is to provide tools to study
discrete hedging which uses random intervention times. We exemplify
these possibilities by using the general results to exhibit a ``no bad
days'' strategy and a minimum standard deviation strategy for the
Black--Scholes model.

Weak convergence theory for approximations of stochastic integrals and
solutions to stochastic differential equations is developed in
\citet{Rootzen}, Kurtz and Protter
(\citeyear{KurtzProttera,KurtzProtterb,KurtzProtterc}) and, in particular, an extensive study of the
Euler method for SDEs is provided by \citet{JacodProtter}. This theory
has been used and extended to solve and analyze various aspects of
approximation and hedging error problems in mathematical finance. As
examples we mention \citet{DuffieProtter}, \citet{BertsimasEtAl},
\citet{HayashiMykland}, \citet{TankovVoltchkova},
\citet{BrodenWiktorsson} and \citet{Fukasawa}. A Malliavin calculus
approach to discrete hedging is used in \citet{GobetTemam} and in a
number of papers, which also consider variable but deterministic grids,
by Geiss and coworkers; see \citet{GeissToivola} and the references
therein. The main theoretical tool of \citet{HayashiMykland} is related
to our general result, as discussed further below. The quite
interesting paper \citet{Fukasawa} also studies random grid
approximations, for one-dimensional processes. The setting of
Fukasawas paper is more or less in the middle between our Theorems
\ref{cormixing} and~\ref{thmerror}. The conditions used by Fukasawa
are rather different from ours, and there does not seem to be any simple
relations between his results and ours.

Now a brief overview of the paper. The next section, Section \ref
{sectgeneral}, contains the basic general theorem on multidimensional
convergence for stochastic integrals with respect to local
multidimensional Brownian semimartingales, and the tools to check
conditions and extend the result. In Section~\ref{sectapproximation} we
give the explicit result for random grid approximations of stochastic
integrals based on the solution of a multidimensional SDE. Section \ref
{sectoptimization} investigates ways to design and optimize
approximation errors, and in Section~\ref{secthedging} this is applied
to discrete financial hedging.

%s2 #&#
\section{General results}
\label{sectgeneral}

This section contains two main results. The first one gives a means to
establish multidimensional convergence of the distribution of
stochastic integrals with more and more rapidly varying integrands, and
the second one shows how convergence of integrals with simple
integrands can be extended to more general integrands. In addition,
Lemma~\ref{lemmascheme} provides tools to check the assumptions of the
theorems. Our main aim is the error in approximations of stochastic
integrals, but the results may in fact also have more general use.

Let $\Omega= C(\R_+, \R^d)$ be the space of continuous $\R^d$-valued
functions defined on $\R_+$, define $B_t= \{B^i_t\}_{1 \leq i \leq d}$
by $B_t(\omega) = \omega(t)$, let $\P$ be the probability measure which
makes $B$ a Brownian motion starting at $0$ and let $\F_t$ be the
completion of the $\sigma$-algebra generated by $\{B_s; 0 \leq s \leq t
\}$. Further write $\F$ for the smallest $\sigma$-algebra which contains
all the $\F_t$. Until further notice is given all random variables we
consider are defined on the filtered probability space $(\Omega, (\F
_t), \F, \P)$. Weak convergence will be for random variables (or
``processes'') with values in $C([0, T], \R^K)$, the space of continuous
$K$-dimensional functions defined on the time interval $[0,T]$, and
with respect to the uniform metric. Usually the dimension $K$ of the
processes will be clear from the context, and then we, for brevity,
write $C[0,T]$, instead of $ C([0, T], \R^K)$, and just write
$\Rightarrow$ for weak convergence.

Weak convergence is stable (or ``R\'enyi-stable'') if it holds on any
subset of $\F$, and the convergence is mixing (or ``R\'enyi-mixing'')
if, in addition, the limit is the same on any subset. In the present
setting this is specified by the definition which follows below. To
appreciate part (ii) of the definition, recall that convergence in
distribution often is written as $X^n \Rightarrow X$, but that in this
notation $X$ is not a random variable defined on some probability
space. It is just a convenient notation for the limiting distribution
of $X_n$. However, one can, of course, construct a random variable with
this distribution, to give $X$ a life of its own.
%
%de2.1 #&#
%
\begin{definition}
\label{defmixing}
(i) Let $(X^n)_{n \geq1}$ be a sequence of random
variables defined on the same probability space $(\Omega, \F, \P)$ and
with values in $C[0, T]$. Then $(X^n)_{n \geq1}$ \textit{converges
stably} if $\E[U f(X^n)]$ converges for any bounded continuous
function $f\dvtx C[0, T] \to\R$ and any bounded measurable random variable
$U$ defined on $(\Omega, \F, \P)$. If, in addition,
%
%e1 #&#
%
\begin{equation}
\label{eqmixing}
\lim_n \E[U f(X^n)] = \E[U] \lim_n\E[f(X_n)],
\end{equation}
then the convergence is \textit{mixing}.

(ii) If $(X^n)_{n \geq1}$ converges stably, then it is always possible
to enlarge the probability space and construct a new random variable
$X$ on the enlarged probability space such that $\lim_n \E[U f(X^n)] =
\E[Uf(X)]$ for all bounded random variables~$U$; see
\citet{AldousEagleson}. Thus, with this construction we can write
stable convergence as $X^n \Rightarrow_{s} X$. If the convergence, in
addition, is mixing, then $X$ is independent of $\F$, and we write $X^n
\Rightarrow_{m} X$.
\end{definition}

It is straightforward to see that to establish stable or mixing
convergence it is enough to prove convergence of $\E[U f(X^n)]$ for
strictly positive $U$ with $EU=1$. Further, see
\citet{AldousEagleson},
$X^n \Rightarrow_{s} X$ \textit{if and only if} $(Y^n, X^n) \Rightarrow
(Y, X^*)$ for any sequence of random variables $Y^n \rightarrow_p Y$
which converges in probability \textit{if and only if} $X^n \Rightarrow X$
with respect to $\P(\cdot| A)$ for any set $A$ with $\P(A)> 0$. (In
the middle statement, convergence is with respect to the product
topology.) Finally, if stability (or mixing) holds with respect to a
sigma-algebra $\F$ and the sigma-algebra $\F'$ is independent of $\F$,
then it also holds with respect to the sigma-algebra generated by $\F$
and $\F'$.

Let $X = (X_j, j= 1, \ldots, d)$ be a continuous $d$-dimensional
Brownian semimartingale defined on the space $(\Omega, (\F_t), \F,
\P)$ by
%
%e2 #&#
%
\begin{equation}
\label{eqsemimartingale}
X_j(t) = \sum_{k=1}^d \int_0^t G_{j,k}(s) \,dB_k(s) + \int_0^t a_j(s)\,ds
\end{equation}
with $G_{j,k}$ and $a_{j}$ adapted, and with $\int_0^T G_{j,k}^2 \,ds<
\infty$ and $\int_0^T a_j^2 \,ds < \infty$ \textit{a.s.} for all $j,k$.
Further\vspace*{1pt} let $\{H_{i,j}^n\} = \{H_{i,j}^n; 1 \leq i,j \leq d\}$ be a
$d\times d$-dimensional array of $\F_t$-adapted processes such that
$\int_0^T (H_{i,j}^n)^2\,dt < \infty$ \textit{a.s.} for each $i, j$,
and write
%
%e3 #&#
%
\begin{eqnarray}
\label{eqarraydef}
\{H_{i,j}^n \cdot X_j \} &=& \{H_{i,j}^n \cdot X_j; 1 \leq i,j \leq d\}
\nonumber\\[-8pt]\\[-8pt]
&=&
\biggl\{\int_0^t H_{i,j}^n(s) \,dX_j(s);
1 \leq i,j \leq d\biggr\}_{0 \leq t \leq T}. \nonumber
\end{eqnarray}
Thus $\{H_{i,j}^n \cdot X_j\}$ takes values in $C([0, T], \R^{d\times
d})$. In the following we let $\to_p$ denote convergence in probability
and take ``positive'' to mean the same as ``nonnegative.''

The form of the second condition, equation \eqref{eqconvofsquarevar} of
the following theorem requires some explanation. For simplicity of
exposition suppressing the index $k$, it says that $\int_0^t H_{i,
j}^nG_{j} H_{l, m}^n G_{m} \,ds$ converges in probability to some
absolutely continuous limit, which we temporarily write as $\int_0^t
C_{(i,j), (l,m)}\,ds$. Since limits of positive variable are positive,
we further assume that for each $t$ and $\omega$ the array $\{C_{(i,j),
(l,m)}(t)\}$ is ``positive definite,'' that is, equivalently, that it
can be obtained as the covariances of some $d \times d$ array of random
variables. The diagonal elements $C_{(i,j), (i,j)}(t)$ of the array are
obtained from the limits of $\int_0^t (H_{i, j}^n)^2G_{j}^2 \,ds$ and
hence it is natural to write them as $C_{(i,j), (i,j)}(t) = (H_{i,
j})^2G_{j}^2$. Further,\vspace*{1pt} taking positive square roots we
may then more generally write $C_{(i,j), (l,m)}(t) = H_{i, j}G_{j}
H_{l, m} G_{m} \rho _{(i,j), (l,m)}$. The array $\{\rho_{(i,j),
(l,m)}\}$ then is the ``correlation array'' corresponding to the
covariances $\{C_{(i,j), (l,m)}(t)\}$. This gives the formulation
\eqref{eqconvofsquarevar}. (If some $ G_{j}$ is zero, we just set the
corresponding $H_{i, j}$'s and off-diagonal elements of $\rho$ to zero,
and the diagonal elements to 1.)

Further, it is possible to find a ``root'' of $\{\rho_{(i,j),
(l,m)}(t)\}$, that is, an array $\{\sigma_{(i,j), (l,m)}(t)\}$ such
that $\rho_{(i,j), (l,m)}(t) = \sum_{1\leq r,s \leq d} \sigma_{(i,j),
(r,s)} \sigma_{(r,s),(l,m)}$. This can be seen by reordering the index
set $\{(i,j); {1\leq r,s \leq d}\}$, linearly, say lexicographically,
making the corresponding reordering of $\{\rho_{(i,j), (i,j)}\}$ into a
matrix which then is positive definite, finding a root of this matrix,
and then making the identification back to the array ordering.
%
%th2.2 #&#
%
\begin{theorem}
\label{cormixing}
Suppose that $\{H_{i,j}^n\}$ satisfies
%
%e4 #&#
%
\begin{equation}
\label{equniform}
\sup_{0 \leq t \leq T}\biggl|\int_0^t H_{i, j}^n \,ds \biggr| \rightarrow_p 0,
\qquad n \to\infty, 1 \leq i,j\leq d,
\end{equation}
and that for $k=1, \ldots, d$
%
%e5 #&#
%
\begin{equation}
\label{eqconvofsquarevar}
\int_0^t H_{i, j}^nG_{j,k} H_{l, m}^n G_{m,k} \,ds \rightarrow_p \int_0^t
H_{i, j}G_{j,k} H_{l, m} G_{m,k}\rho^{k}_{(i,j), (l, m)} \,ds
\end{equation}
as $n \to\infty$, for $i,j,l,m=1, \ldots, d$, and for some correlation
array processes $\rho^{k} =(\rho^{k}_{(i,j), (l, m)}; k= 1,
\ldots, d)$ and processes $\{H_{i, j}\dvtx 1 \leq i, j \leq d \}$ such that all
$H_{i, j}G_j$ are positive. Let $\sigma^{k}(t)$ be an arbitrary root of
$\rho^{k}(t)$; see the discussion just before the theorem.
Then, for $X$ given by \eqref{eqsemimartingale},
%
%e6 #&#
%
\begin{equation}
\label{eqgenconv}
\{H_{i,j}^n \cdot X_j\} \Rightarrow_s \Biggl\{\sum^{d}_{r,s,k=1}
H_{i,j}G_{j,k} \sigma_{(i,j), (r, s)}^{k} \cdot W_{r,s,k}\Biggr\}
\end{equation}
as $n \to\infty$, where $W =(W_{r,s,k}; 1 \leq r,s,k \leq d)$ is a
$d\times d \times d$-dimensional Brownian motion which is independent
of $\F$.
\end{theorem}

This result simplifies in the special case when $X$ is just a Brownian
motion $B$; see the following corollary. The corollary is close to
Theorem A.1 of \citet{HayashiMykland}. Differences are that the
corollary makes the basic condition \eqref{equniform} explicit, gives a
more detailed description of the limit distribution and has the more
powerful conclusion of stable convergence.

In Theorem~\ref{cormixing} we, for simplicity of notation, considered a
quadratic array $\{H_{i,j}^n \cdot X_j\dvtx 1 \leq i,j \leq d \}$. This
does\vspace*{1pt} not involve any loss of generality, but still, for later use in
the proof of Theorem~\ref{cormixing}, it is convenient to formulate the
corollary for a rectangular array.
%
%co2.3 #&#
%
\begin{corollary}
\label{thmmixing}
Suppose that \eqref{equniform} is satisfied for $i=1, \ldots, d_1,
j=1, \ldots, d_2$
and that
%
%e7 #&#
%
\begin{equation}
\label{eqtau1}
\int_0^t H_{i, k}^n H_{j, k}^n \,ds \rightarrow_p \int_0^t H_{i, k} H_{j,
k}\rho_{i,j}^{k} \,ds,\qquad n \to\infty,
\end{equation}
as $n \to\infty$, for some correlation matrix processes $\rho^{k}=
\sigma^{k}(\sigma^{k})'$, where $i,j=1, \ldots, d_1, k=1, \ldots,
d_2$, and positive processes $\{H_{i, k}\dvtx i= 1, \ldots, d_1, k=1,
\ldots, d_2\}$, and for $0 \leq t \leq T$.
Then
%
%e8 #&#
%
\begin{equation}
\label{eqstableconvergence}
\{H_{i,k}^n \cdot B_k \} \Rightarrow_s \Biggl\{\sum^{d_1}_{j=1} H_{i, k}
\sigma^{k}_{i,j} \cdot W_{j,k} \Biggr\}
\end{equation}
as $n \rightarrow\infty$, where $W =\{W_{j,k}\dvtx j= 1, \ldots, d_1,
k=1, \ldots, d_2\}$ is a Brownian motion which is independent of $\F$.
\end{corollary}

The following lemma plays an important role in the proofs.
%
%le2.4 #&#
%
\begin{lemma}
\label{lemmaconvtozero}
Suppose that $\eta(t)$ and $H^n(t)$ are real-valued random processes
with $\int_0^S \eta(t)^2 \,dt < \infty$ a.s. and with
$\limsup
_{n \to\infty} \int_0^S H^n(t)^2 \,dt < \infty$ a.s. for some
positive constant $S \leq\infty$. Suppose further that
\[
\sup_{0 \leq t \leq S}\biggl|\int_0^tH^n \,ds \biggr| \rightarrow_p 0,
\qquad n\to\infty.
\]
Then
%
%e9 #&#
%
\begin{equation}
\label{equniformzero}
\sup_{0 \leq t \leq S}\biggl|\int_0^tH^n\eta \,ds \biggr| \rightarrow_p 0,
\qquad n \to\infty.
\end{equation}
\end{lemma}
\begin{pf}%{Proof of Lemma~\ref{lemmaconvtozero}}
Suppose first that there exists a sequence $\{\eta_k\}$ of processes
such that
\begin{eqnarray*}
\label{eqetaappr}
\int_0^S \bigl(\eta(t) - \eta_k(t)\bigr)^2 \,dt &\rightarrow_p& 0
\qquad\mbox{as }
k \to\infty, \\
\sup_{0 \leq t \leq S}
\biggl|\int_0^t H^n\eta_k(s) \,ds\biggr| &\rightarrow_p& 0
\qquad\mbox{as } n \to\infty\qquad \mbox{for each }k.
\end{eqnarray*}
Then, by the Cauchy--Schwarz inequality,
\begin{eqnarray*}
\limsup_n \sup_{0 \leq t \leq S} \biggl|\int_0^t H^n\eta \,ds\biggr| &\leq&
\limsup_n
\sup_{0 \leq t \leq S} \biggl|\int_0^t H^n\eta_k \,ds\biggr| \\
&&{} + \limsup_n \sup_{0 \leq t \leq S}\biggl|\int_0^t H^n(\eta- \eta_k)
\,ds\biggr| \\
&\leq& 0 + \sqrt{\limsup_n \int_0^S (H^n)^2 \,dt} \sqrt{\int_0^S
(\eta-
\eta_k)^2 \,dt},
\end{eqnarray*}
which tends to 0 as $k \to\infty$, so that \eqref{equniformzero} holds.

Thus the lemma follows if there exist a sequence $\{\eta_k\}$ which
satisfies the two requirements above.

Now, for each $k$ there exists a continuous process $\tilde{\eta}_k $,
measurable in $t$ and $\omega$, such that $\P(\int_0^S(\eta(t) -
\tilde
{\eta}_k(t))^2 \,dt > 1/k) \leq1/k$. Briefly, to see this note that if
$\eta(t)$ is approximated by convolving it with a sequence of
``approximate $\delta$-functions,'' for example, with a sequence of
centered normal densities with variance parameters tending to $0$, then
the convolutions are measurable in $t$ and $\omega$ and for almost all
$\omega$ converge to $\eta[\cdot, \omega)$ in $L_2[0,S]$. The existence
of the sequence $\tilde{\eta}_k$ follows at once from this, since
convergence \textit{a.s.} implies convergence in probability.

Next, with $1_A$ denoting the indicator function of a set $A$, for
$\tilde{\eta}_{k,m}(t) = \sum_{i=0}^{[mS]} \tilde{\eta}_{k}(iS/m) 1_{\{t
\in[ iS/m, (i+1)S/m)\}}$ it follows that
\[
\int_0^S\bigl(\tilde{\eta}_k(t) - \tilde{\eta}_{k, m}(t)\bigr)^2 \,dt
\rightarrow
_{\mathrm{a.s.}} 0 \qquad\mbox{as } m \to\infty
\]
and thus, choosing $m_k$ suitably, $\eta_k = \sum_{i=0}^{[m_kS]}
\tilde
{\eta}_k( iS/m_k) 1_{\{t \in[ iS/m_k, (i+1)S/m_k)\}}$ satisfies the first
one of the two relations above. Furthermore, the second one is easily
seen to hold for $\eta_k$ of this form.
\end{pf}
\begin{pf*}{Proof of Theorem~\ref{cormixing} and Corollary
\ref{thmmixing}}
We do this in reverse order, and first prove Corollary~\ref{thmmixing}.
For simplicity of notation we only prove the corollary for a
two-dimensional Brownian motion, that is, for the case $d=2$. The general
case is the same.

By\vspace*{1pt} Rootz{\'e}n [(\citeyear{Rootzen}), Theorem 1.2], each marginal process $\{H_{i,
j}^{n}\cdot B_j(t), 0 \leq t \leq T\}$ is tight $C([0,T],\R)$, and then
also the entire $d\times d$-dimensional sequence $\{H_{i, j}^{n}\cdot
B_j (t), 0 \leq t \leq T, 1\leq i,j \leq d\}$ is tight\vspace*{1pt} $C([0,T],\R
^{d\times d})$, so only stable finite-dimensional convergence remains
to be proved. We prove this in two steps, where the first one follows
along the lines of \citet{Rootzen} and the second step uses the
Cram\'er--Wold device. A final third step uses Corollary~\ref{thmmixing} to
prove Theorem~\ref{cormixing}.\vspace*{8pt}

\textit{Step} 1: Let $\{\psi_i^n; i=1, 2\}$ be adapted processes such
that, for $i=1, 2$,
%
%e10 #&#
%
\begin{equation}
\label{eqpsiunifzero}
\sup_{0 \leq t \leq T} \biggl|\int_0^t \psi_i^n \,ds\biggr| \rightarrow_p 0
\end{equation}
and such that
%
%e11 #&#
%
\begin{equation}
\label{eqpsisquareconv}
\int_0^t (\psi_i^n)^2 \,ds \to_p \int_0^t (\psi_i)^2 \,ds
\end{equation}
for some $\psi_1, \psi_2 > 0, 0 \leq t \leq T$. To make inverses
well defined, we, without loss of generality, can assume that the $\psi
_i^n(t)$ are defined also for $t >T$, and such that equations \eqref
{eqpsiunifzero} and \eqref{eqpsisquareconv} hold with $T$ replaced by
$S$ for any $S>0$, and with \mbox{$\psi_i(t)=1$} for $t > T$ and $i=1, 2$.
This does not involve the result to be proved nor the assumptions, and
hence can be done without loss of generality.

Let $C[0, \infty) = C([0, \infty), \R)$ be the space of continuous real
valued functions defined on $[0, \infty)$ and endowed with the topology
of uniform convergence on compact sets; see \citet{Whitt}. Let the
random variable $U>0$ satisfy $\E U = 1$, and assume the functional
$f\dvtx C[0, \infty) \to\R$ is bounded and continuous. Further, set $\tau_n(t)
= \int_0^t (\psi_1^n)^2 \,ds + \int_0^t (\psi_2^n)^2 \,ds$, let $\tau
(t) =
\lim_{n \to\infty} \tau_n(t) = \int_0^t (\psi_1)^2 \,ds + \int_0^t
(\psi
_2)^2 \,ds$ and define $\tau_n^{-1}$ by $\tau_n^{-1} (t) = \inf\{s\dvtx\tau
_n(s) > t \}$. Additionally let $\tilde{W}$ be a one-dimensional Brownian
motion which is independent of $\F$. We first prove that
%
%e12 #&#
%
\begin{equation}
\label{equf}
\E U f \biggl(\int_0^{\tau_n^{-1} (\cdot)} \psi_1^n \,dB_1 + \int_0^{\tau
_n^{-1} (\cdot)} \psi_2^n \,dB_2\biggr)
\to\E f(\tilde{W}(\cdot)),
\end{equation}
for each such $U$, so that $\int_0^{\tau_n^{-1} (\cdot)} \psi_1^n \,dB_1
+ \int_0^{\tau_n^{-1} (\cdot)} \psi_2^{n} \,dB_2 \Rightarrow_m
\tilde
{W}$, on $C[0, \infty)$.

Now, define a new probability measure $\Q$ by $d\Q/d\P= U$, and write
$\E_\Q$ for expectation taken with respect to $\Q$. Then, by Girsanov's
theorem [\citet{rogerswilliams2001},\vadjust{\goodbreak} Theorem IV 38.5] there exists an
adapted square integrable process $c = (c_1, c_2)$ such that $(\tilde
{B}(t) = (B_1(t) - \int_0^t c_1(s)\,ds, B_2(t) - \int_0^t c_2(s)\,ds)$ is a
Brownian motion under $\Q$.

Hence,
%
%e13 #&#
%
\begin{eqnarray}
\label{eqgirsanov}
&&\E U f \biggl(\int_0^{\tau_n^{-1} (\cdot)} \psi_1^n \,dB_1 + \int_0^{\tau
_n^{-1} (\cdot)} \psi_2^n \,dB_2\biggr) \nonumber\\
&&\qquad= \E_{\Q} f \biggl(\int_0^{\tau_n^{-1} (\cdot)} \psi_1^n \,d\tilde{B}_1
+ \int
_0^{\tau_n^{-1} (\cdot)} \psi_2^n \,d\tilde{B}_2 \\
&&\hspace*{26.6pt}\qquad\quad{} + \int_0^{\tau_n^{-1} (\cdot)} \psi_1^n c_1 \,ds +\int
_0^{\tau
_n^{-1} (\cdot)} \psi_2^n c_2 \,ds\biggr).\nonumber
\end{eqnarray}
Under $\Q$ the process $ \int_0^{\tau_n^{-1} (\cdot)} \psi_1^n
\,d\tilde
{B}_1 + \int_0^{\tau_n^{-1} (\cdot)} \psi_2^n \,d\tilde{B}_2$ has the
same distribution as $\tilde{W}$ [\citet{rogerswilliams2001}, Theorem IV
34.1]. Further, by Lemma~\ref{lemmaconvtozero}, we have that $\int_0^t
\psi_1^n c_1 \,ds +\int_0^t \psi_2^n c_2 \,ds \to_p 0$ in $C[0, S]$, for
any fixed $S$. Since $f$ is bounded and continuous on $C[0, \infty)$,
these two facts prove \eqref{equf}, and hence mixing convergence on
$C[0, \infty)$.

It thus follows from $\tau_n \rightarrow_p \tau$ that $(\tau_n,
\int
_0^{\tau_n^{-1} (\cdot)} \psi_1^n \,dB_1 + \int_0^{\tau_n^{-1}
(\cdot)}
\psi_2^n \,dB_2)) \Rightarrow_s (\tau, \tilde{W})$, and hence, by
composing $\tau_n^{-1}$ with $\tau_n$ [cf. \citet{Billingsley}, page
145], that
%
%e14 #&#
%
\begin{equation}
\label{eqconvtowtilde}
\int_0^t \psi_1^n \,dB_1 + \int_0^t \psi_2^n \,dB_2 \Rightarrow_s
\tilde
{W}(\tau(t))
\end{equation}
in $C[0, \infty)$, and hence, in particular, in $C[0,T]$.\vspace*{8pt}

\textit{Step} 2: Finite-dimensional stable convergence now follows by
standard but notationally complicated Cram\'er--Wold arguments. To
lessen complications we here only consider two basic cases, and leave
the general argument to the reader. Thus, first, let $\psi_i^n(s) =
b_i1_{\{0 \leq s \leq t_i\}}H_{1,i}^n(s)$ for $i=1, 2$, with $0<t_1,
t_2 \leq T$. Equation \eqref{eqtau1} implies that
\[
\tau_n(t) \rightarrow_p \tau(t) = b_1^2\int_0^{t \wedge t_1}(H_{1,1})^2
\,ds + b_2^2\int_0^{t \wedge t_2}(H_{1,2})^2 \,ds
\]
so that by \eqref{eqconvtowtilde},
\begin{eqnarray*}
&&
b_1 \int_0^{t \wedge t_1} H_{1,1}^n \,dB_1 + b_2\int_0^{t \wedge
t_2}H_{1,2}^n \,dB_2 \\
&&\qquad\Rightarrow_s \tilde{W}\biggl( b_1^2\int_0^{t \wedge
t_1}(H_{1,1})^2 \,ds + b_2^2\int_0^{t \wedge t_2}(H_{1,2})^2 \,ds\biggr).
\end{eqnarray*}

Now,\vspace*{1pt} using elementary properties of Brownian motion together with
Rogers and Williams [(\citeyear{rogerswilliams2001}), Theorem IV 34.1]
we have that $\tilde{W}( b_1^2\int_0^{t \wedge t_1}(H_{1,1})^2 \,ds +
b_2^2\int_0^{t \wedge t_2}(H_{1,2})^2 \,ds)$ has the same distribution,\vadjust{\goodbreak}
and the same dependency with any $\F $-measurable variable, as
\[
b_1\int_0^{t \wedge t_1} H_{1,1} \,dW_{1,1} + b_{2}\int_0^{t \wedge t_2}
H_{1,2} \,dW_{1,2}
\]
for independent Brownian motions $W_{1,1}, W_{1,2}$, so that we by
\eqref{eqconvtowtilde} have established that $b_1 \int_0^{t \wedge t_1}
H_{1,1}^n \,dB_1 + b_2\int_0^{t \wedge t_2}H_{1,2}^n \,dB_2 \Rightarrow_s
b_1 \int_0^{t \wedge t_1} H_{1,1} \,dW_{1,1} + b_2\int_0^{t \wedge
t_2}H_{1,2} \,dW_{1,2}$, for any real numbers $b_1, b_2$. In particular
stable two-dimen\-sional convergence of $(H_{1,1}^n\cdot B_1(t_1),
H_{1,2}^n\cdot B_2(t_2))$ to $(\int_0^{t_1} H_{1,1} \,dW_{1,1},\break \int
_0^{t_2} H_{1,2} \,dW_{1,2})$ follows by Cram\'er--Wold.

If we instead take $\psi_1^n = b_1I_{\{0 \leq s \leq t_1\}}H_{1,1}^n(s)
+ b_2I_{\{0 \leq s \leq t_2\}}H_{2,1}^n(s)$ and $\psi_2^n = 0$ then, by
\eqref{eqtau1},
\begin{eqnarray*}
\tau_n(t) &\rightarrow_p& \tau(t)\\
&=& b_1^2\int_0^{t \wedge t_1}(H_{1,1})^2
\,ds +2b_1b_2 \int_0^{t \wedge t_1 \wedge t_2}H_{1,1}H_{2,1}\rho_{1,
2}^1\,ds\\
&&{}+ b_2^2\int_0^{t \wedge t_2}(H_{2,1})^2 \,ds.
\end{eqnarray*}
Furthermore, similarly as before and recalling that the matrix $\sigma
^1$ is a root of the correlation matrix $\rho^1$, it can be seen that
then $\tilde{W}(\tau(\cdot))$ has the same distribution, and the same
dependency with any $\F$-measurable variable, as
\begin{eqnarray*}
&&b_1 \biggl(\int_0^{t \wedge t_1}H_{1,1} \sigma^{1}_{1,1}\,dW_{1,1} + \int_0^{t
\wedge t_1}H_{1,1}\sigma_{1,2}^1\,dW_{2,1}\biggr) \\
&&\qquad{}
+ b_2 \biggl(\int_0^{t \wedge t_2}H_{2,1}\sigma_{2,1}^1\,dW_{1,1} + \int_0^{t
\wedge t_2} H_{2,1}\sigma^{1}_{2,2}\,dW_{2,1} \biggr).
\end{eqnarray*}
Reasoning as above we get that
\begin{eqnarray*}
&& b_1 \int_0^{t_1} H_{1,1}^n \,dB_1 + b_2\int_0^{t_2}H_{2,1}^n \,dB_1
 \\
&&\qquad\Rightarrow_s b_1\int_0^{t_1} H_{1,1}\sigma^{1}_{1,1} \,dW_{1,1} + b_1\int
_0^{t_1} H_{1,1}\sigma^{1}_{1,2} \,dW_{2,1} \\
&&\qquad\quad{} + b_2 \int_0^{t_2} H_{2,1}\sigma^{1}_{2,1} \,dW_{1,1} + b_2
\int
_0^{t_2} H_{2,1}\sigma^{1}_{2,2} \,dW_{2,1}
\end{eqnarray*}
for independent Brownian motions $W_{1,1}, W_{2,1}$. Since $b_1$ and
$b_2$ are arbitrary, this proves stable two-dimensional convergence of
$(H_{1,1}^n\cdot B_1(t_1), H_{2,1}^n\cdot B_1(t_2))$. A general proof
of Corollary~\ref{thmmixing} is only notationally more complicated.

We next use Corollary~\ref{thmmixing} to obtain the conclusion of
Theorem~\ref{cormixing}.\vspace*{6pt}

\textit{Step} 3:
By Lemma~\ref{lemmaconvtozero}, if $H_{i,j}^n$ satisfies \eqref
{equniform}, then $\sup_{0 \leq t \leq T} |\int_0^t H_{i,j}^n a_i\,ds|
\to
_p 0$, for all $i, j$, and hence the general\vadjust{\goodbreak} result follows if we can
prove that the result of the theorem holds for the case when all $a_i$
are identically zero. Thus, to find the limit of $\{H_{i,j}^n \cdot
X_j\}$ one only has to consider
\[
\Biggl\{\sum_{k=1}^d H_{i,j}^nG_{j,k} \cdot B_k\Biggr\}.
\]
Again by Lemma~\ref{lemmaconvtozero}, if $H_{i,j}^n$ satisfies \eqref
{equniform}, then
%
%e15 #&#
%
\begin{equation}
\sup_{0 \leq t \leq T}\biggl|\int_0^t H_{i, j}^{n}G_{j,k}\,ds\biggr| \to_{p} 0.
\end{equation}
Now, making the definition $ H_{(i,j), k}^n:= H_{i,j}^nG_{j,k}$ and
replacing the index $i$ in \eqref{eqstableconvergence} by the
``multiindex'' $(i,j)$, convergence of the array $ \{H_{i,j}^nG_{j,k}
\cdot B_k\}
$ follows from Corollary~\ref{thmmixing} with $d_1 =d^2, d_2=d$. The
result \eqref{eqgenconv} then follows by summing over $k$ and writing
$W_{l,m,k}$ for $W_{(l,m),k}$.
\end{pf*}

We now change to a more general setup, from Brownian semimartingales to
general processes $(H^n,X^n)$ which are defined on filtered probability
spaces $\Psi^n = (\Omega^n,\F^n,\P^n,(\F^n_t)_{0\leq t < \infty})$.
Here $\F^n$ is a $\P^n$-complete $\sigma$-algebra and $(\F
^n_t)_{0\leq
t < \infty}$ is a filtration which satisfies the usual hypotheses (but
which is not necessarily generated by a Brownian motion). The following
definition is key to our goal. We give it for vector valued processes.
The definition for matrix valued processes is analogous.
%
%de2.5 #&#
%
\begin{definition}
\label{defgoodness}
Let $(X_n)_{n\geq1}$ be a sequence of continuous $\R^d$-valued
semimartingales defined on $\Psi^n$, $n\geq1$ and assume that $X^n
\Rightarrow X$. The sequence $X^n$ is \textit{good} if for any sequence
of $\R^{d\times d}$-valued adapted c\`adl\`ag stochastic processes
$(H^n)_{n \geq1}$ defined on $\Psi^n$ such that $(H^n,X^n)
\Rightarrow
(H,X)$, there exists a filtration $(\mathcal{G}_t)$ such that $X$ is a
semimartingale and $H$ is an adapted c\`adl\`ag process, and $\{
H^n_{i,j} \cdot X^n_j \} \Rightarrow\{ H_{i,j}\cdot X_j \}$.
\end{definition}

The following criterion is sufficient for goodness; see, for example,
Theorem~2.2 in \citet{KurtzProttera}.
%
%de2.6 #&#
%
\begin{definition}
\label{ucv}
A sequence of continuous $\R^d$-valued semimartingales $(X^{n})_{n\geq
1}$ is said to have \textit{uniformly controlled variations} (UCV) if
for each $n \geq1$, there exist decompositions $X^{n} = M^{n} + A^{n}$
such that
\[
\sup_n \E^n\biggl\{ [M^{n},M^{n}]_{T} + \int_0^{T} |dA_s^{n}|
\biggr\}
< \infty.
\]
\end{definition}

The next theorem combined with Theorem~\ref{thmmixing} will give the
asymptotic distributions of approximation errors for stochastic
integrals. If, in addition to the conditions of the theorem,\vadjust{\goodbreak} $f$ is
bounded, then the result follows from Theorem 3.5 in \citet
{KurtzProtterb}. However, in the present setting the result holds also
without the boundedness condition, and it is further possible to give
a quite simple proof. In the theorem, $0=\tau_0^n < \tau_1^n < \cdots<
\infty$ are $\{ \F_t\}$-stopping times, and $\eta_n$ is defined by
$\eta
_n(t) = \tau_k^n$, $\tau_k^n \leq t < \tau_{k+1}^n$.
%
%th2.7 #&#
%
\begin{theorem}
\label{thmmain}
Let $Y$ be a continuous $\R^d$-valued $\{ \F_t\}$-semimartingale on
$[0, T]$, and suppose that $f = (f_1, \ldots, f_d)$ is continuously
differentiable. Assume that $\eta_n(t)$ tends to the identity in
probability for $t \in[0, T]$, and let $\{ \lambda_n \}$ be a positive
sequence converging to infinity. Further, set
\begin{eqnarray*}
U^n &=& \lambda_n\int\bigl(f(Y) - f(Y \circ\eta_n)\bigr)\,dY \\
:\!&=& \lambda_n \sum_{i=1}^d\int\bigl(f_i(Y) - f_i(Y \circ\eta_n)\bigr)\,dY_i
\end{eqnarray*}
and define
%
%e16 #&#
%
\begin{equation}
\label{Zn}
Z_{ij}^n(t) = \lambda_n \int_0^t \bigl(Y_i(s) - Y_i\circ\eta_n(s)\bigr)\,dY_j(s).
\end{equation}
Suppose that $(Z^n)_{n \geq1}$ is good, and that $(Z^n, Y) \Rightarrow
(Z, Y)$. Then $U^n \Rightarrow U$ on $[0, T]$,
where
\[
U = \sum_{i,j=1}^d \int\frac{\partial f_j(Y)}{\partial y_i}\,dZ_{ij}.
\]
\end{theorem}

Since $\eta_n$ is nondecreasing, pointwise convergence in probability
in $[0, T]$, as assumed in the theorem, is equivalent to uniform
convergence in probability in $[0, T]$. Below we will use this without
further comment.
\begin{pf*}{Proof of Theorem~\ref{thmmain}}
For simplicity of exposition, we assume that $d=1$.
By the continuous mapping theorem we have that $(Z^n,Y,Y) \Rightarrow
(Z,Y,Y)$. Since $Y$ is continuous, and $\eta_n$ converges uniformly in
probability to the unity, this in turn can be seen to imply that
$(Z^n,Y \circ\eta_n,Y) \Rightarrow(Z,Y,Y)$, for example, by using the
Skorokhod translation of convergence in distribution to convergence
\textit{a.s.}

We now define
\[
g(x,y) = \frac{f(x)- f(y)}{x - y},
\]
where we make the continuous choice $g(x,x) = f'(x)$ when the
denominator vanishes. The function $g$ is uniformly continuous on $[0,
T]^2$, so the continuous mapping theorem gives that $(Z^n,g(Y,Y\circ
\eta_n) ) \Rightarrow(Z,f'(Y))$. Now,
\[
U^n = \lambda_n\int\bigl(f(Y) - f(Y \circ\eta_n)\bigr)\,dY = \int g(Y,Y \circ
\eta_n)\,dZ^n.
\]
But since $(Z^n)_{n\geq1}$ is good, we have that
\[
\int g(Y,Y\circ\eta_n) \,dZ^{n} \Rightarrow\int f'(Y)\,dZ,
\]
which proves the theorem for $d=1$.
\end{pf*}

The next lemma provides a tool for verification of criteria like \eqref
{equniform} and \eqref{eqtau1}. In the lemma we specialize to stopping
times (cf. the \hyperref[intro]{Introduction}) defined recursively by
$\tau_0^n = 0$ and
%
%e17 #&#
%
\begin{equation}
\label{eqtaukn}
\tau_{k+1}^n = \biggl(\tau_k^n + \frac{1}{n\theta(\tau_k^n)}\biggr)
\wedge T
\end{equation}
for some adapted stochastic process $\theta$. As before, let
%
%e18 #&#
%
\begin{equation}
\label{eqeta}
\eta_n(t) = \tau_k^n, \qquad\tau_k^n \leq t < \tau_{k+1}^n\qquad
\mbox{for } k=1, 2, \ldots
\end{equation}
and write $E_p = \E\int_0^1 B(s)^p \,ds = \int_0^1 s^{p/2} \E B(1)^p
\,ds =
\E B(1)^p /(p/2+1)$ so that $E_1 = \E\int_0^1 B(s) \,ds = 0$ and $E_2 =
\E
\int_0^1 B(s)^2 \,ds = 1/2$.

In the lemma we will assume that the function $a(t); t \in[0, T]$ is
\textit{locally bounded}, that is, that to any $\varepsilon> 0$ there
exists a localizing stopping time $\nu= \nu_\varepsilon$ such that $a(t
\wedge\nu); t \in[0, T]$ is bounded, and such that $\P(\nu< T) <
\varepsilon$. In particular, if $a$ is continuous on $[0, T]$, then $a$ is
locally bounded.
%
%le2.8 #&#
%
\begin{lemma}
\label{lemmascheme}
Assume that $a$ and $\theta$ are adapted processes such that $a$ is
locally bounded, $\theta$ is strictly positive and $a(t)/\theta
(t)^{p/2}$ is a.s. Riemann integrable over $[0, T]$, and let $\tau
_k^n$ and $\eta_n$ be defined by \eqref{eqtaukn} and \eqref{eqeta}. Set
%
%e19 #&#
%
\begin{equation}
\psi_n (t) = n^{p/2}\sum_{k=0}^{\infty}a(\tau_k^n)\bigl(B(t) - B(\tau
_k^n)\bigr)^p 1_{\{\tau_k^n \leq t < \tau_{k+1}^n \}}.
\end{equation}
Further assume that $\eta_n$ tends to the identity in probability. Then
%
%e20 #&#
%
\begin{equation}
\label{eqpsiapprox}
\sup_{0 \leq t \leq T} \biggl|\int_0^t \psi_n(s) \,ds - E_p \int_0^t \frac
{a(s)}{\theta(s)^{p/2}}\,ds\biggr| \rightarrow_p 0
\end{equation}
as $n \rightarrow\infty$, for $p = 1, 2$.
\end{lemma}
\begin{pf}%{Proof of Lemma~\ref{lemmascheme}}
If we prove the lemma under the additional restriction that $a$ is
bounded, then it follows in general, since it then holds for $a(t)$
replaced by $a(t \wedge\nu)$ for any localizing stopping time $\nu$,
and this in turn implies that \eqref{eqpsiapprox} holds with
probability greater than $1-\varepsilon$, for arbitrary $\varepsilon$. Thus
we assume in the rest of this proof that $a$ is uniformly bounded, so
that in particular the expectations exist.

To ease notation we below sometimes will write $\tau_k$ instead of
$\tau
_k^n$ and define $\bar{\F}_k = \F_{\tau_k}$. Clearly
\begin{eqnarray*}
n^{p/2}\E\biggl\{\int_{\tau_k}^{\tau_{k+1}}a(\tau_k)\bigl(B(t) - B(\tau
_k)\bigr)^p \,dt\Big|
\bar{\F}_k
\biggr\} &=& n^{p/2}a(\tau_k) \int_0^{1/n\theta(\tau_k)} \E B(t)^p \,dt \\
&=& E_p \frac{a(\tau_k)}{n \theta(\tau_k)^{p/2 +1}}.
\end{eqnarray*}
Recalling the definition of $\eta_k$,
\[
\sum_{k'=1}^{k-1} E_p \frac{a(\tau_{k'})}{n \theta(\tau_{k'})^{p/2
+1}} =
E_p \int_0^{\tau_k} \frac{a \circ\eta_n (s)}{\theta\circ\eta_n
(s)^{p/2}}\,ds
\]
and hence
\[
X_k:= \int_0^{\tau_k} \psi_n \,ds - E_p \int_0^{\tau_k} \frac{a
\circ
\eta_n (s)}{\theta\circ\eta_n (s)^{p/2}} \,ds
\]
is a martingale with index set $\mathbb Z_+$.

In the following we show that $\sum_k \E((X_{k+1} - X_k)^2 |
\bar
{\F}_k ) \rightarrow0$. By the functional central limit theorem
for martingales [see, e.g., \citet{Rootzen1983}, Theorem 3.5] this in
turn implies that
%
%e21 #&#
%
\begin{equation}
\label{eqmartingaletozero}
{\max_k}|X_k| = {\max_{k }} \biggl|\int_0^{\tau_k^n}\psi_n \,ds - E_p \int
_0^{\tau
_k^n} \frac{a \circ\eta_n (s)}{\theta\circ\eta_n (s)^{p/2}} \,ds\biggr|
\rightarrow_p 0
\end{equation}
as $n \rightarrow\infty$.
Using the Cauchy--Schwarz inequality in the second step, elementary
properties of Brownian motion in the third and that $(\tau_{k+1} -
\tau
_k) = 1/(n\theta(\tau_{k}))$ in the fourth step, we have that
\begin{eqnarray*}
&&
\sum_k \E[(X_{k+1} - X_k)^2 | \bar{\F}_k ] \\
&&\qquad\leq
\sum_k \E
\biggl[\biggl(\int_{\tau_k}^{\tau_{k+1}} \psi_n \,dt \biggr)^2 \Big|
\bar{\F}_k
\biggr] \\
% &=& n^p \sum_k a(\tau_k)^2 \E[ \int_{\tau_k}^{\tau_{k+1}} (B(t)
%- B(\tau_k))^{2p} \,dt | \bar{\F_k}] \\
&&\qquad\leq n^p \sum_k a(\tau_k)^2 (\tau_{k+1} - \tau_k) \int_{\tau
_k}^{\tau
_{k+1}} \E\bigl[\bigl(B(t) - B(\tau_k)\bigr)^{2p}| \bar{\F}_k\bigr] \,dt
\\
&&\qquad= \frac{E_{2p}}{p+1}n^p \sum_k a(\tau_k)^2(\tau_{k+1} - \tau
_k)^{p+2} \\
&&\qquad\leq \frac{E_{2p}}{p+1} \max_k\biggl(\frac{a(\tau_k)}{n\theta(\tau
_k)^{p/2+1}}\biggr) \sum_{k\dvtx \tau_k < T} \biggl(\frac{a(\tau_k)}{n\theta(\tau
_k)^{p/2+1}}\biggr).
\end{eqnarray*}
It follows from the Riemannn integrability of $a/\theta^{p/2}$ that in
the last expression above the first factor tends to $0$ and that the
second tends to $\int_0^T a(s)/\theta(s)^{p/2}\,ds$, so that the product
tends to zero. This completes the proof of \eqref{eqmartingaletozero}.

The assumption that $a$ is bounded and straightforward computation show
that $\E\int_0^T\psi_n^2 \,ds$ is bounded in $n$, and since furthermore
$\max_k \{ \tau_{k+1}^n - \tau_k^n \} \to_p 0$, we can apply the
Cauchy--Schwarz inequality, to see that
\[
\max_k \sup_{\tau_k^n \leq t < \tau_{k+1}^n} \biggl|\int_{\tau_k^n}^t
\psi_n
\,ds\biggr| \leq\biggl( \max_k \{\tau_{k+1}^n - \tau_k^n \}\int_0^T\psi_n^2
\,ds\biggr)^{1/2} \rightarrow_p 0
\]
for $n \rightarrow\infty$. Together with \eqref{eqmartingaletozero}
this shows that
%
%e22 #&#
%
\begin{equation}
\label{eqintervalzero}
\sup_{0 \leq t \leq T} \biggl|\int_0^{t}\psi_n \,ds - E_p \int_0^{t} \frac{a
\circ\eta_n (s)}{\theta\circ\eta_n (s)^{p/2}} \,ds\biggr| \rightarrow_p 0.
\end{equation}
By assumption $a/\theta^{p/2}$ is Riemann integrable, and hence
%
%e23 #&#
%
\begin{equation}
\label{eqriemannconv}
\sup_{0 \leq t \leq T} \biggl|\int_0^{t} \frac{a(s)}{\theta(s)^{p/2}} \,ds -
\int_0^{t} \frac{a \circ\eta_n (s)}{\theta\circ\eta_n (s)^{p/2}} \,ds\biggr|
\rightarrow_{\mathrm{a.s.}} 0.
\end{equation}
The triangle inequality together with \eqref{eqintervalzero} and
\eqref
{eqriemannconv} completes the proof of the lemma.
\end{pf}

%s3 #&#
\section{Approximation of stochastic integrals}
\label{sectapproximation}

We now use the results from the previous section to find the explicit
form of the asymptotic distribution of the sum of the errors in
approximating $d$ stochastic integrals where the integrands are
functions of the solution to a $d$-dimensional SDE and where the
integrators are the same solutions to the SDE. The following condition
is used in the theorem.
%
%co3.1 #&#
%
\begin{condition}
\label{condsde}
Let the measurable functions $\alpha(\cdot)\dvtx \R^d \rightarrow\R^d,
\beta(\cdot)\dvtx \R^d \rightarrow\R^{d \times d}$ satisfy
\[
|\alpha(x)| + |\beta(x)| \leq C(1 + |x|),
\]
where $x \in\R^d$ for some constant $C$ and
\[
|\alpha(x) - \alpha(y)| + |\beta(x) - \beta(y)| \leq D|x-y|,
\]
where $x,y \in\R^d$ for some constant $D$.
\end{condition}

This condition ensures that the SDE has an unique continuous solution.
Further, we will need the following lemma, which is given as Lemma 2.5
in \citet{Rootzen1983}.
%
%le3.2 #&#
%
\begin{lemma}
\label{lemmarootzen}
Suppose $\{Z^{n}\}_{n\geq1}$ is a sequence of positive discrete time
stochastic processes, adapted to their respective filtrations $\{\F
^{n}\}_{n\geq1}$ and that $\tau_{n}$ is a stopping time with respect
to $\F^{n}$ for each $n$. Then
\[
\sum_{j=1}^{\tau_{n}}\E(Z^{n}_{j}|\F^{n}_{j-1}) \rightarrow_{p} 0\vadjust{\goodbreak}
\]
implies that
\[
\sum_{j=1}^{\tau_{n}}Z^{n}_{j} \rightarrow_{p} 0.
\]
\end{lemma}
%
%th3.3 #&#
%
\begin{theorem}
\label{thmerror}
Let $Y$ be the solution of the SDE
%
%e24 #&#
%
\begin{equation}
\label{eqdiff}
dY(t) = \alpha(Y(t)) \,dt + \beta(Y(t))\,dB(t),
\end{equation}
where $B$ is a $d$-dimensional Brownian motion, $\alpha, \beta$ satisfy
Condition~\ref{condsde} and $Y(0)$ is independent of $B$ and
satisfies $\E Y(0)^2 < \infty$. Then the error in the Euler-type
approximation scheme defined by
\begin{eqnarray*}
U^n(t) &=& n^{1/2}\int_0^t \bigl(f(Y(u)) - f\bigl(Y \circ\eta
_n(u)\bigr)
\bigr)\,dY(u) \\
:\!&=& n^{1/2} \sum_{i=1}^d \int_0^t \bigl(f_i(Y(u)) - f_i\bigl(Y \circ
\eta
_n(u)\bigr)\bigr)\,dY_i(u),
\end{eqnarray*}
where $f=(f_1, \ldots, f_d)$ is continuously differentiable and the
grid is given by \eqref{eqtaukn} with $\sup_{t\in[0,T]}\theta
(t)<\infty
$ \textit{a.s.} and $1/\theta$ \textit{a.s.} Riemann integrable, satisfies
\[
U^n \Rightarrow\sum_{r, k = 1}^d \int_0^t \Delta_{r,k}(u)\,dW_{r,k}(u)
\]
on $[0, T]$, where
\[
\Delta_{r,k}(t) = \frac{\sum_{i, j = 1}^d ({\partial
f_j}/{\partial
y_i})(Y(t))\beta_{i,r}(Y(t)) \beta_{j,k}(Y(t))}{\sqrt{2\theta(t)}},
\]
and $W$ is an $d\times d$-dimensional Brownian motion, independent of
$B$. In particular,
\[
\sup_{0 \leq t \leq T} |U^n(t)| \Rightarrow\sup_{0 \leq t \leq T}\Biggl|
\sum
_{r,k = 1}^d \int_0^t \Delta_{r,k}(u)\,dW_{r,k}(u)\Biggr|.
\]
\end{theorem}
\begin{pf}%{Proof of Theorem~\ref{thmerror}}
For the convenience of the reader we begin by recalling that $\{X^n\}
_{n\geq1}$ is $O_p(a_n)$ for some sequence $a_n$ if
\[
\lim_{c \to\infty}\limsup_{n \to\infty}\P[ |X^n/a_n|\geq c ]=0
\]
or, equivalently, if $\{X_n/a_n\}_{n\geq1}$ is \textit{tight}.
We first assume that the coefficients $\alpha$ and $\beta$ are
uniformly bounded, and prove that the result holds under this extra
assumption. The general result for unbounded coefficients then follows
by an easy localization argument which is given at the end of the
proof. We again write $\bar{\F}_{v}$ instead of $\F_{\tau_{v}}$ and
often suppress the explicit dependence on $n$ and, for example, write
$\tau_{v}$ instead of $\tau_{v}^n$.\vadjust{\goodbreak}

Since $1/\theta$ is Riemann integrable, and hence pathwise bounded
\textit{a.s.}, and $\sup_{t\in[0,T]}\theta(t)<\infty$ \textit {a.s.},
it\vspace*{1pt} follows that $\eta_n$ tends to $t$ uniformly
\textit{a.s.} By Theorem~5.2.1 in \citet{Oksendal} there exists a
unique $t$-continuous solution $Y$ to equation \eqref{eqdiff}.

The first part of the proof consists of proving that
\[
\{Z^{n}_{i,j}\} = \biggl\{ \sqrt{n}\int_0^t \bigl(Y_i(s) - Y_i \circ\eta
_n(s)\bigr)\,dY_j(s) \biggr\}
\]
converges jointly with $Y$. We do this by showing that the conditions
of Theorem~\ref{cormixing} are satisfied for the choices $H_{i,j}^n =
\sqrt{n} (Y_i - Y_i\circ\eta_n )$ and $G_{j, k} = \beta_{j, k}$.

The bounded variation part of $Y_i - Y_i \circ\eta_n$ can be seen to
give contributions which are $O_p(1/n)$, and thus, using the triangle
inequality and writing $1_v(s) = 1_{\{\tau_v \leq s < \tau_{v+1}\}}$,
it can be seen that \eqref{equniform} follows if we show that
%
%e25 #&#
%
\begin{equation}
\label{eqsdeuniformtozero}
\sqrt{n}\sup_{t\in[0,T]}\biggl|\int^{t}_{0}\sum_{v}\int_{\tau_{v}}^s
1_v(s) \beta_{i,j}(u)\,dB_{j}(u)\,ds\biggr| \rightarrow_{p}0
\end{equation}
for $1\leq i,j \leq d$.

Now,
%
%e26 #&#
%
\begin{eqnarray}
\label{eqerrorparts}
&&
\sqrt{n} \int^{t}_{0} \sum_{v}\int_{\tau_v}^{s}1_{v}(s)\beta
_{i,j}(u)
\,dB_{j}(u)\,ds \nonumber\\
&&\qquad= \sqrt{n} \int^{t}_{0}\sum_{v}\int_{\tau_v}^{s}1_{v}(s)\bigl(\beta
_{i,j}(u) - \beta_{i,j}(\tau_{v})\bigr)\,dB_{i}(u)\,ds \\
&&\qquad\quad{} +\sqrt{n} \int
^{t}_{0}\sum_{v}1_{v}(s)\beta_{i,j}(\tau_{v})\bigl(B_{i}(u) - B_{i}(\tau
_{v})\bigr)\,ds. \nonumber
\end{eqnarray}
The last term tends to zero in probability by Lemma~\ref{lemmascheme}
with $p=1$, since Riemann integrability of $1/\sqrt{\theta}$ follows
from Riemann integrability of $1/\theta$.

We next show that also the first term on the right-hand side is negligible.
Let $C$ denote a generic deterministic constant whose value may change
from one appearance to the next. Since $\tau_{v+1}$ is measurable with
respect to $\bar{\F}_v$ it follows from Condition~\ref{condsde},
It\^o's isometry, and the assumption that the constants in \eqref{eqdiff}
are bounded that
%
%e27 #&#
%
\begin{eqnarray}
\label{eqsqrvar}
\E\biggl[\int_{\tau_v}^{s}\bigl(\beta_{i,j}(u) - \beta
_{i,j}(\tau_{v})\bigr)^2\,du\Big|\bar{\F}_v\biggr] &\leq& C \int_{\tau_v}^{s} \E
[|Y(u) -
Y(\tau_v)|^2 | \bar{\F}_v]\,du \nonumber\\
&\leq& C \int_{\tau_v}^{s}(u-\tau_v) \,du \\
&\leq& C(\tau_{v+1} - \tau_v)^2.
\nonumber
\end{eqnarray}
Define
\[
\Delta_v(t) = \sqrt{n} \int^{t\wedge\tau_{v+1}}_{\tau_v}\int
_{\tau
_v}^{s\wedge\tau_{v+1}}\bigl(\beta_{i,j}(u) - \beta_{i,j}(\tau_{v})\bigr)\,dB_{i}(u)\,ds,
\]
so that the first term on the right-hand side of \eqref{eqerrorparts}
equals $\sum_v \Delta_v(t)$.
Using Doob's inequality together with the Cauchy--Schwarz inequality in
the second step and \eqref{eqsqrvar} in the third step we have that
\begin{eqnarray*}
&&
\E\Bigl[\sup_{\tau_v \leq t < \tau_{v+1}} |\Delta_v(t)| \big| \bar{\F}_v\Bigr]\\
&&\qquad\leq
\sqrt{n}(\tau_{v+1} - \tau_v) \E\biggl[\sup_{\tau_v \leq s < \tau_{v+1}}\biggl|
\int_{\tau_v}^{s}\bigl(\beta_{i,j}(u) - \beta_{i,j}(\tau
_{v})\bigr)\,dB_{i}(u)\biggr| \Big|
\bar{\F}_v\biggr] \\
&&\qquad\leq C \sqrt{n} (\tau_{v+1} - \tau_v) \E\biggl[ \int_{\tau_v}^{\tau
_{v+1}}\bigl(\beta_{i,j}(u) - \beta_{i,j}(\tau_{v})\bigr)^2 \,du \Big| \bar{\F
}_v\biggr]^{1/2}\\
&&\qquad\leq C \sqrt{n} (\tau_{v+1} - \tau_v)^2.
\end{eqnarray*}
Thus, by the definition \eqref{eqtaukn},
\begin{eqnarray*}
\sum_v \E\Bigl[\sup_{\tau_v \leq t < \tau_{v+1}} |\Delta_v(t)| \big| \bar
{\F
}_v\Bigr] &\leq& C \sqrt{n} \sum_v (\tau_{v+1} - \tau_v)^2 \\
&\leq& C \sqrt{n}\frac{1}{n} T\sup_{0 \leq t \leq T } \frac
{1}{\theta
(t)} \to_{\mathrm{a.s.}} 0.
\end{eqnarray*}
According to Lemma~\ref{lemmarootzen} it follows that $\sum_v \sup
_{\tau
_v \leq t < \tau_{v+1}} |\Delta_v(t)| \to_p 0$. Hence,
\[
\sup_{0 \leq t \leq T} \biggl|\sum_v \Delta_v(t)\biggr| \leq\sum_v \sup_{\tau_v
\leq t < \tau_{v+1}} |\Delta_v(t)| \to_p 0,
\]
which completes the proof that the first term in the right-hand side of
\eqref{eqerrorparts} tends uniformly to zero in probability.

Completely similar, but more complex computation show that for any
indexes $i,j,k,l,m$, and using Lemma~\ref{lemmascheme} with $p=2$ for
$j=m$ and computations similar to (but simpler than) the proof of Lemma
\ref{lemmascheme} for $j \neq m$,
\begin{eqnarray*}
&&
n\int^{t}_{0} \sum_{v}\int_{\tau_{v}}^{s}1_{v}(s)\beta_{i,j}
(u)\,dB_{j}(u)
\int_{\tau_{v}}^{s}1_{v}(s)\beta_{l,m}(z)\,dB_{m}(z)\beta
_{j,k}(s)\beta
_{m,k}(s)\,ds \\
&&\qquad= n\int^{t}_{0} \sum_{v}\beta_{i,j}(\tau_{v})\beta_{j,k}(\tau
_{v})\beta
_{l,m}(\tau_{v})\beta_{m,k}(\tau_{v})\\
&&\hspace*{67.5pt}{}\times\bigl(B_j(s) - B_j( \tau_{v} )\bigr)\bigl(B_m(s) - B_m(\tau_{v} ) \bigr) 1_v(s) \,ds
+ o_p(1)\\
&&\qquad\rightarrow_p \frac{1}{2}\int^{t}_{0} \beta_{i,j}(s)\beta_{j,k}(s)
\beta_{l,m}(s)\beta_{m,k}(s)/\theta(s) \delta_{j,m} \,ds,
\end{eqnarray*}
where\vspace*{1pt} $\delta_{j,m}$ is 1 if $j=m$ and zero otherwise. Recalling that
$G_{j, k} = \beta_{j, k}$, and as before approximating\vspace*{1pt} $H_{i,j}^n(s) =
\sqrt{n} (Y_i(s) - Y_i\circ\eta_n (s))$ by $\sum_{k=1}^d \sum
_{v}\int
_{\tau_{v}}^{s}1_{v}(s)\*\beta_{i,k}(u)\,dB_{k}(u)$ it follows that
condition \eqref{eqconvofsquarevar} of Theorem~\ref{cormixing} holds as
%
%e28 #&#
%
\begin{equation}
\label{eqsdeqv}
\int^{t}_{0}H^{n}_{i,j}G_{j,k}H^{n}_{l,m}G_{m,k}\,ds \rightarrow_{p}
\frac
{1}{2}\sum^{d}_{r=1}\int^{t}_{0}\beta_{i,r}\beta_{l,r}\beta
_{j,k}\beta
_{m,k}/\theta \,ds.
\end{equation}
Now we recognize that the choice $H_{i,j}G_{j,k}\sigma
^{k}_{(i,j),(r,s)}= \delta_{r,s} \beta_{j,k}\beta_{i,r}/\sqrt
{2\theta}
$ satisfies equation \eqref{eqsdeqv}. Hence,
%
%e29 #&#
%
\begin{equation}
\label{sdestable}
\{Z^{n}_{i,j}\}=\{H_{i,j}^{n} \cdot Y_j\} \Rightarrow_s \Biggl\{ \sum_{r,k =
1}^d \frac{\beta_{j,k}\beta_{i,r}}{\sqrt{2\theta}} \cdot W_{r,k} \Biggr\}.
\end{equation}

Arguments similar to those above show that $\{H_{i,j}^{n} \cdot Y_j\}$
has uniformly controlled variations and hence are good. Stable
convergence implies that the left-hand side of \eqref
{sdestable} converges jointly with $Y$. The first conclusion of the
theorem now follows from Theorem~\ref{thmmain}, for the case when the
coefficients are bounded.

To remove the restriction that the coefficients are bounded, for
general $\alpha_i, \beta_{i,j}$ define coefficients $\alpha_{i}^{c}=
(-c)\vee\alpha_{i} \wedge c$ and $\beta_{i,j}^{c}= (-c)\vee\beta_{i,j} \wedge c$.
Theorem~5.2.1 in \citet{Oksendal} still yields unique $t$-continuous
solution $Y^{c}$ to \eqref{eqdiff} for these functions. Let
$U^{n,c}$ be defined from $\alpha_{i}^{c}, \beta_{i,j}^{c}$ in the same
way as $U^n$ is defined from $\alpha_{i}, \beta_{i,j}$. With obvious
notation, we have already proved that $U^{n,c}\Rightarrow U^{c}$, as $n
\to\infty$ for each fixed $c$. Since $\P(\sup_{t\in[0,T]}|Y^{c}(t) -
Y(t)|>0)\rightarrow0$, as $c \rightarrow\infty$ also
$U^{c}\Rightarrow U$. Further,
\begin{eqnarray*}
&&\limsup_{n}\P\Bigl({\sup_{t\in[0,T]}}|U^{n,c} - U^{n}|> 0\Bigr)\\
&&\qquad \le\P\bigl(\inf\bigl\{t\dvtx\max\{\max\{|\alpha_i(\bar{Y}_t)|\},\max\{|\beta
_{i,j}(\bar{Y}_t)|\}\}\ge c\bigr\} \le T\bigr) \rightarrow0
\end{eqnarray*}
as $c \rightarrow\infty$. Hence, Theorem 3.2 in \citet{Billingsley}
gives that $U^{n} \Rightarrow U$, which proves that the first result of
the theorem holds also for the general case.

The second conclusion follows from from the first by the continuous
mapping theorem, since the supremum mapping is continuous.
\end{pf}

%s4 #&#
\section{Designing the error in approximations of stochastic integrals}
\label{sectoptimization}

In deciding on which approximation scheme to use to compute a
stochastic integral---or, to decide on a hedging strategy---one has to
balance the error with the \textit{number of intervention times} $N =N_n=
\max\{k; \tau_k^{n} < T\}$. In this section we will investigate two
such schemes. The first one could be called the ``no bad days''
strategy, and simply consists in choosing the stopping times $\{ \tau_k
\}$ where the stochastic integral is evaluated---or the times when the
portfolio is rehedged---in such a way that the error is a Wiener
process. In the second strategy we bound the expected number of
evaluation times and minimize the asymptotic standard deviation of the
approximation error under this restriction.

The setting of this section is the following: suppressing the
superscript $n$ the stopping times are given by \eqref{eqtaukn}, that is,
$\tau_0 = 0$ and
%
%e30 #&#
%
\begin{equation}
\label{eqtaukk}
\tau_{k+1} = \biggl(\tau_k + \frac{1}{n\theta(\tau_k)}\biggr)
\wedge T
\end{equation}
with $\theta$ adapted and positive, and the distribution of the
approximation error $\varepsilon(t)$ satisfies
%
%e31 #&#
%
\begin{equation}
\label{eqerror}
\sqrt{n}\varepsilon(t) \Rightarrow\int_0^t \frac{f(s)}{\sqrt{\theta
(s)}} \,dW(t)
\end{equation}
for some adapted process $f(s) \geq0$ and Wiener process $W$ which is
independent of $\theta$ and $f$. Here it should be noted that \eqref
{eqerror} is more general than it looks at first; for example, the
approximation error in Theorem~\ref{thmerror} satisfies this for $f(t)
= \sqrt{\frac{1}{2}\sum_{k,m = 1}^d \Delta_{k,m}^2(t)}$.

It is straightforward to find the asymptotic number of evaluation times.
%
%pr4.1 #&#
%
\begin{proposition}
\label{propinterventiontimes}
Suppose that $\theta$ is Riemann integrable a.s. and that $\inf_{0
\leq
t \leq T}\theta(t) > 0$ a.s. Then
\[
\lim_{n\rightarrow\infty} \frac{N_n}{n} = \int_0^T \theta(t)\,dt\qquad
\mbox{a.s.}
\]
If, in addition, $\E[\sup_{0 \le t \le T}\theta(t)]<\infty$, then
\[
\lim_{n\rightarrow\infty} \E\frac{N_n}{n} = \int_0^T \E\theta(t)\,dt.
\]
\end{proposition}
\begin{pf}%{Proof of Proposition~\ref{propinterventiontimes}}
Suppose first $\theta$ is of the form
%
%e32 #&#
%
\begin{equation}
\label{eqthetastep}
\theta(t) = \sum_{i=0}^k \theta_i 1_{[a_i, a_{i+1})}(t)
\end{equation}
for some random variables $\theta_i>0$ and constants $0=a_0 < a_1 < \cdots<
a_k=T$, and with $1_{[a_i, a_{i+1})}$ the indicator function of the
interval $[a_i, a_{i+1})$. For each $\omega$, it is easily seen that
the number of intervention times in the interval $ [a_i, a_{i+1})$ is
$n\theta_i(a_{i+1}-a_i) + O(1)$, and hence
\[
\frac{N_n}{n} = \sum_{i=0}^k \theta_i(a_{i+1}-a_i) + O\biggl(\frac{1}{n}\biggr) =
\int_0^T \theta(t)\,dt + O\biggl(\frac{1}{n}\biggr) \to\int_0^T \theta(t)\,dt
\]
as $n \to\infty$. If $ \tilde{\theta} \leq\theta$ and $\tilde
{\theta
}$ is of the form \eqref{eqthetastep} then, with obvious notation,
$N_n(\tilde{\theta}) \leq N_n(\theta) + O(1)$, and the corresponding
bound with all the inequalities reversed is also true.

Now, by assumption $\theta$ is Riemann integrable, and hence can be
approximated arbitrarily well from below and above by functions of the
form \eqref{eqthetastep}. This proves the first assertion of the proposition.

Furthermore, $N_n/n \leq T \sup_{0 \leq t \leq T}\theta(t) + 1/n$, and
hence the second assertion follows from the first one by dominated convergence.
\end{pf}

In the rest of this section we assume that we ``are in the asymptotic
regime,'' that is, that $n$ is so large that we, to the degree of
approximation needed, may assume that the limits above can be replaced
by equalities. Thus, below we will assume that
%
%e33 #&#
%
\begin{equation}\label{eqaserror}
\E N = n\int_0^T \E\theta(t)\,dt,\qquad \varepsilon(t) = \frac
{1}{\sqrt{n}}\int_0^t \frac{f(s)}{\sqrt{\theta(s)}} \,dW(t),
\end{equation}
so that in particular $\E\varepsilon(t)^2 = \frac{1}{n}\int_0^t \E
\frac
{f(s)^2}{\theta(s)}\,ds$.\vspace*{8pt}

\textit{The no bad days strategy}: It is at once seen, supposing that
$f^2$ is Riemann integrable, that if we choose $\theta(t) = cf(t)^2$,
for some constant $c$, then
\[
\varepsilon(t) = \frac{1}{\sqrt{cn}}W(t)
\]
and
\[
\E N = cn\int_0^T \E f^2(s)\,ds.
\]
Thus, in a financial setting, with this choice of $\theta$, there are
no ``days'' where the hedging error grows quicker than during other
days, and hence a trader can sleep equally well (or equally badly!)
each night.\vspace*{8pt}

\textit{Minimal standard deviation}: We will now, supposing that $f$ is
Riemann integrable, show that the solution of the optimization problem
\[
\inf_{\{\theta\dvtx \theta\geq0, \mathrm{adapted}\}} \bigl\{\sqrt{\E
\varepsilon
^2(T)}\dvtx \E N \leq nC\bigr\}
\]
is given by $\theta(t)= Cf(t)/(\int_0^T \E f(s)\,ds$). For this choice
\[
\E N = nC,\qquad \varepsilon(t) = \sqrt{\frac{\int_0^T \E
f\,ds}{nC}}\int_0^t \sqrt{f} \,dW.
\]
Thus in particular, for the optimal strategy the standard deviation
is\break
$\sqrt{\E\varepsilon(T)^2} = \int_0^T E f\,ds/ \sqrt{nC}$.

Now, write $\tilde{\theta} = n \theta$. With this notation $\E
\varepsilon
(T)^2 = \E\int_0^T f^2/\tilde{\theta}\,ds$ and the restriction is $\E
\int
_0^T \tilde{\theta} \leq nC$. Applying the Cauchy--Schwarz inequality
twice, it follows that
\begin{eqnarray*}
\biggl(\E\int_0^T f \,ds\biggr)^2 &\leq& \biggl( \E\sqrt{\int_0^T f^2/\tilde
{\theta
}\,ds} \sqrt{ \int_0^T \tilde{\theta}\,ds}\biggr)^2 \\
&\leq& \E\biggl(\int_0^T f^2/\tilde{\theta}\,ds\biggr) \E\biggl(\int_0^T \tilde
{\theta
}\,ds\biggr)
\end{eqnarray*}
and hence
\[
\E\varepsilon(t)^2 \geq\frac{(\E\int_0^T f \,ds)^2}{nC}.
\]
However, above we have seen that $\theta= Cf/(\int_0^T \E f \,ds)$
achieves this bound, and hence is the optimal choice.

%{\em Minimal expected absolute error:} This is the same optimization
%problem as above, but with $\E\varepsilon(T)^2$ replaced %by
%ds}.
%We have not been able to solve this problem in general. However,
%reasoning as before it can be seen that the minimal %standard
%deviation strategy also leads to smaller expected absolute error than
%the "no bad days" strategy. In fact the choice %$c=1/(n\int_0^T \E
%f^2(s)\,ds)$ in the no bad days strategy leads to $\E N =C$ and $\E|
%However, for the minimum variance strategy the expected value of the
%absolute value of the error %is
% \sqrt{\frac{2}{\pi C}} \E\int_0^T f \,ds \leq\sqrt{\frac{2 }{\pi C}}\E
% which establishes this claim.

%s5 #&#
\section{Application to hedging}
\label{secthedging}

An important application of the results in the previous section is to
hedging of financial derivatives. Here we treat the simplest
Black--Scholes model and only give a brief comment on more complicated
problems. The limit distribution of the Black--Scholes hedging error
for equidistant deterministic grids has been studied, for example, in
\citet{BertsimasEtAl} and \citet{HayashiMykland}. [We have not been able
to follow the proof of Theorem 1.b in \citet{BertsimasEtAl};
specifically, we could not understand the use of Lemma 5.1 from \citet
{DuffieProtter}.]

We distinguish between \textit{complete} and \textit{incomplete}
financial markets. In complete markets, all derivatives can be
replicated (hedged) perfectly by trading in a self-financing way in the
underlying and a money market account. The approximation error
distribution we analyze is here the total hedging error. In an
incomplete market, an investor who hedges a contract will still choose
a hedging portfolio which is, in some sense, optimal for her purposes.
In this case, the error we obtain is relative to this optimal hedging
portfolio. We give now an application of the results in the previous
section to hedging in the complete Black--Scholes market.

We assume that a stock $S$ follows the Black--Scholes model. In other
words, we model the stock as a geometric Brownian motion, which has the dynamics
\[
dS(t) = \mu S(t)\,dt + \sigma S(t)\,dB(t)
\]
for $\mu, \sigma> 0$, where $B$ is a Brownian motion, and $S(0) =
s>0$. Further, we have a risk-free money market account with dynamics
\[
dR(t) = rR(t)\,dt
\]
for $r>0$, where $R(0)=1$. It is well known that the price of a
so-called \textit{call option} with payoff $\max(S(T) - K, 0)$ at the
deterministic terminal\vadjust{\goodbreak} time $T$, for some strike price $K$, is at time $t$
\[
\Pi(t) = \Phi(d_+)S(t) - Ke^{-(T-t)}\Phi(d_-),
\]
where $\Phi$ denotes the standard normal cumulative distribution
function and
\[
d_{\pm}(t) = \frac{\log({S(t)}/{K}) + (r \pm{\sigma
^2}/{2})(T-t)}{\sigma\sqrt{T-t}}.
\]
Now, if we set
\[
Y(t) = \pmatrix{
S(t) \cr
R(t)}
\]
and $f = (\Phi(d_+), -\Phi(d_-)Ke^{-rT})$, we get that
\[
\Pi(t) = \int f(Y(t))\,dY(t)
\]
gives the self-financing price process of the call option. This is of
the form considered in Theorem~\ref{thmerror}, with $d=2$ and $\beta
_{1,1}(t) = \sigma S(t)$, and all other $\beta$-s equal to zero. Thus,
using the stopping times \eqref{eqtaukn}, Theorem~\ref{thmerror} gives
that the hedging error satisfies
\begin{eqnarray*}
\sqrt{n} \bigl(\Pi(t) - \Pi\circ\eta_n(t)\bigr) &\Rightarrow& \int_0^t
\frac
{df_1}{dx_1}(s) \sigma^2 S(s)^2/\sqrt{2\theta(s)} \,dW(s) \\
&=& \int_0^t \frac{\phi(d_+(t))\sigma S(s)}{\sqrt{2\theta
(s)(T-s)}} \,dW(s)
\end{eqnarray*}
with $\phi(t) = d\Phi(t)/dt$ the standard normal density function.

Consider now an investor who hedges a call option, but who only adjusts
her hedge at some stopping times $\{\tau_k\}_{k \geq1}$ of her own
choosing. If she wants to have a ``uniform'' increase of the error and
make it approximately a Brownian motion, she should use the ``no bad
days'' strategy from the previous section. This would mean that she
would use the stopping times \eqref{eqtaukk} with $\theta(t) = c\phi
(d_+(t))^2 \sigma^2 S(t)^2/(2(T-t))$. However, this leads to a (purely)
technical difficulty: $\theta(t)$ tends to $0$ as $t \to T$ if $S(T)
\in\R\setminus K$ and to $\infty$ if $S(T)=K$. This means that the
assumption of \textit{a.s.} Riemann integrability of $1/\theta$ is not
satisfied on $[0,T]$, nor is the assumption that $\sup_{t\in
[0,T]}\theta
(t)<\infty$ on $[0,T]$. A theoretical (and in fact also practical)
solution is to instead only evaluate the hedging strategy up to a
constant time $V<T$, with $V$ close to $T$. Theorem~\ref{thmerror}
gives that the hedging error up until $V$ for large $cn$ then
approximately is distributed as $W(t)/\sqrt{cn}$.

Alternatively, the minimum standard deviation strategy and the same
reasoning as above lead to choosing
%
%e34 #&#
%
\begin{equation}
\theta(t) = \frac{C\phi(d_+(t)) \sigma S(t)}{\sqrt{2(T-t)}}\Big/\biggl(n\int_0^V
\E\biggl[\frac{\phi(d_+(s)) \sigma S(s)}{\sqrt{2(T-s)}}\biggr]\,ds\biggr),
\end{equation}
where $C$ is the expected number of evaluation times. This yields the
approximate distribution
%
%e35 #&#
%
\begin{equation}
\sqrt{\int_0^V \E\biggl[\frac{\phi(d_+(s)) \sigma S(s)}{C\sqrt
{2(T-s)}}\biggr]\,ds}\int_0^{t} \sqrt{\frac{\phi(d_+(s)) \sigma S(s)}{\sqrt
{2(T-s)}}} \,dW(s)
\end{equation}
for the hedging error, for $n$ large.

It is now completely straightforward to add one or more stocks to the
portfolio and, using, for example, that $\int_0^t f_1 /\sqrt{\theta}
\,dW_1 + \int_0^t f_2/\sqrt{\theta} \,dW_2$ has the same distribution as
$\int_0^t \sqrt{f_1^2 + f_1^2} /\sqrt{\theta} \,dW$, to find the optimal
stopping times and the resulting error when the hedges for all of the
stocks are adjusted at the same time points. This is how portfolio
hedging is done in practice. We leave these calculations to the reader.

An alternative and equally interesting application of our results is
to the field of portfolio optimization. For example, in managing a
large equity portfolio a tracking error arises due to that it is
expensive, or otherwise infeasible, to rebalance the portfolio back to
its optimal state too frequently. Since the optimal portfolio to be
held by the investor is always known, we are exactly in the setting of
the present paper. Here, too, we leave the calculations to the reader.

\section*{Acknowledgment}

We thank an anonymous referee for a very helpful reading which has led
to very substantial improvement.

%suskaldyti doi

% imsref loaded by lrinkeviciute, 2012-05-22 07:59:34
% imsref loaded by lrinkeviciute, 2012-05-22 08:02:06

\printaddresses

\end{document}